\documentclass[12pt,twoside]{amsart}

\newtheorem{Theorem}{Theorem}

\newtheorem{Conjecture}[Theorem]{Conejecture}
\newtheorem{Lemma}[Theorem]{Lemma}

\theoremstyle{remark}

\newtheorem{Definition}{Definition}

\numberwithin{equation}{section}

\begin{document}

\newbox\Adr
\setbox\Adr\vbox{
\vskip18pt
\centerline{Department of Studies in Mathematics}
\centerline{University of Mysore, Manasagangotri}
\centerline{Mysore -- 570 006, Karnataka, India}
\centerline{E-mail: \tt padma\_vathamma@yahoo.com, chandru\_alur@yahoo.com,}
\centerline{\tt raghu\_maths@yahoo.co.in}
\vskip18pt
\centerline{*Institut Girard Desargues,
Universit\'e Claude Bernard Lyon-I}
\centerline{21, avenue Claude Bernard, F-69622 Villeurbanne Cedex, France}
\centerline{E-mail: \tt kratt@igd.univ-lyon1.fr}
\centerline{WWW: \footnotesize\tt http://igd.univ-lyon1.fr/\~{}kratt}
}

\title[Analytic proof of a partition identity]{Analytic
proof of the partition identity
$A_{5,3,3}(n) = B^0_{5,3,3}(n)$}
\author[Padmavathamma, Chandrashekara, Raghavendra,
 Krattenthaler]{Padmavathamma, Chandrashekara B.M., Raghavendra R.,
and C.~Krattenthaler*$^\dagger$\\[18pt]\box\Adr}
\address{Department of Studies in Mathematics, University of Mysore,
Manasagangotri, Mysore -- 570 006, Karnataka, India}
\address{Institut Girard Desargues,
Universit\'e Claude Bernard Lyon-I,
21, avenue Claude Bernard,
F-69622 Villeurbanne Cedex, France}
\thanks{$^\dagger$Research partially supported by
EC's IHRP Programme, grant HPRN-CT-2001-00272,
``Algebraic Combinatorics in Europe"}



\begin{abstract}
In this paper we give an analytic proof
of the identity $A_{5,3,3}(n) =B^0_{5,3,3}(n)$,
where $A_{5,3,3}(n)$ counts the number of partitions of $n$ subject to
certain restrictions on their parts,
and $B^0_{5,3,3}(n)$ counts the number of
partitions of $n$ subject to certain other restrictions on their parts,
both too long to be stated in the abstract. Our proof establishes
actually a refinement of that partition identity.
The original identity was first
discovered by the first author
jointly with M.~Ruby Salestina and S.~R.~Sudarshan in
[``A new theorem on partitions," Proc.\ Int.\ Conference on Special
Functions, IMSC, Chennai, India, September 23--27, 2002; to appear],
where it was also given a combinatorial proof,
thus responding a question of Andrews.
\end{abstract}
\maketitle

\section{Introduction} \label{sec:1}

For an even integer $\lambda$, let
$A_{\lambda ,k,a}(n)$ denote the number of partitions of  $n$
such that

\begin{itemize}
\item  no part $\not\equiv 0$ (mod $\lambda+1$) may
be repeated, and
\item  no part is $\equiv0,\pm \left(a -
\frac{\lambda}{2}\right)(\lambda +1)$ (mod
$(2k-\lambda+1)(\lambda+1)$).
\end{itemize}

For an odd integer $\lambda$, let $A_{\lambda ,k,a}(n)$ denote
the number of partitions of $n$ such that

\begin{itemize}
\item no part  $ \not \equiv 0$ \big(mod
$\frac{\lambda+1}{2}$\big) may be repeated,
\item no part is $\equiv
\lambda +1$ (mod $2\lambda +2$), and
\item  no part is $\equiv
0,\pm(2a-\lambda)\left(\frac{\lambda +1}{2}\right)$ (mod
$(2k-\lambda+1)(\lambda+1)$).
\end{itemize}

\medskip
Let $B_{\lambda ,k,a}(n)$
denote the number of partitions of $n$ of the form $b_1 + \cdots
+ b_{s}$ with $b_{i} \geq b_{i+1}$, such that

\begin{itemize}
\item no part $\not \equiv 0$ (mod
$\lambda+1)$ is repeated,
\item $b_{i} - b_{i+k-1} \geq \lambda+1$,
with strict inequality if $b_i$ is a multiple of $\lambda+1$, and
\item $\sum_{i=j}^{\lambda-j+1}f_{i} \leq a-j$ for $1 \leq j \leq
\frac{\lambda+1}{2}\;$ and $f_1+\dots+f_{\lambda+1} \leq a-1$, where
$f_{j}$ is the number of appearances of $j$ in the
partition.
\end{itemize}

\medskip
In 1969, Andrews \cite{R1} proved the following theorem.

\begin{Theorem}[{\protect\cite[Th.~2]{R1}}] \label{thm:1}
 \ If $\lambda$, $k $, and $a$ are
 positive integers with $  \frac{\lambda}{2} \leq a
\leq k$ and $ k \geq 2\lambda - 1  $, then, for every positive
integer $n$, we have $$A_{\lambda ,k,a}(n)=B_{\lambda ,k,a}(n).$$
\end{Theorem}

Schur's theorem \cite{R6} addresses the case $  \lambda = k = a = 2$.
Hence, it is {\it not\/} a particular case of Theorem~\ref{thm:1} as $ k \geq
2\lambda  - 1  $ is not satisfied. Motivated by this observation,
Andrews \cite{R1}
first conjectured, and later proved in \cite{R2}, that Theorem~\ref{thm:1}
is still true for $  k \geq \lambda.$

In the same paper \cite{R2},
Andrews raised the following question: Is it possible to modify the
conditions on the partitions enumerated by $B_{\lambda ,k,a}(n)$ so
that values of $k < \lambda$
would be admissible? In fact, Schur \cite{R6} had proved that
$A_{3,2,2}(n) = B^0_{3,2,2}(n)$,
where $B^0_{3,2,2}(n)$ denotes the number of partitions enumerated by
$B_{3,2,2}(n)$ with the
added condition that no parts are $\equiv 2$ (mod 4).

This led  Andrews \cite{R2} to state the following conjecture.

\begin{Conjecture}\label{conj}
There holds the identity $A_{4,3,3}(n) = B^0_{4,3,3}(n) $ for all
positive  integers $n$,  where $B^0_{4,3,3}(n)$ denotes the
number of partitions of $ n $ enumerated by $B_{4,3,3}(n)$ with
the added restrictions:
\begin{align*}
f_{5j+2} + f_{5j+3} &\leq 1\quad \text{for}\quad j\ge0,\\
f_{5j+4} + f_{5j+6} &\leq 1\quad \text{for}\quad j\ge0,\\
f_{5j-1} + f_{5j} + f_{5j+5} + f_{5j+6} &\leq 3 \quad \text{for}\quad j\ge1,
\end{align*}
where, as before,
$f_{j}$ denotes the number of appearances of $j$ in the
partition.
\end{Conjecture}

In the year 1994, Andrews et al.\ \cite{R3} gave an analytic proof of
the above conjecture. The first
author and Ruby Salestina, M. gave a combinatorial
proof in \cite{R4}. In \cite{R5},
these two authors and Sudarshan, S.R.
first conjectured, and then proved combinatorially, the following result,
which is analogous to Conjecture~\ref{conj}.

\begin{Theorem}\label{thm:2}
There holds the identity $A_{5,3,3}(n) = B^0_{5,3,3}(n) $ for all
positive  integers $n$,  where $B^0_{5,3,3}(n)$ denotes the
number of partitions of $ n $ enumerated by $B_{5,3,3}(n)$ with
the added restrictions:
\begin{align*}
 f_{6j+3} &= 0\quad \text{for}\quad j \geq 0,\\
 f_{6j+2} + f_{6j+4} &\leq 1\quad \text{for}\quad j\ge0,\\
f_{6j+5} + f_{6j+7}&\leq 1\quad \text{for}\quad j\ge0,\\
f_{6j-1} + f_{6j} + f_{6j+6} + f_{6j+7}&\leq 3 \quad \text{for}\quad j\ge1.
\end{align*}
\end{Theorem}

     The object of this paper is to give an analytic proof of the
partition identity stated in
Theorem~\ref{thm:2}. Actually, we are going to prove a new refinement
of that partition identity, which we state in the next section.
The method of our proof in Section~\ref{sec:3}
is similar to that of Andrews et al.\ in \cite{R3}.

\section{A refinement of the partition identity in
Theorem~\ref{thm:2}} \label{sec:2}

Before being able to state the announced refinement of
Theorem~\ref{thm:2}, we need to make two definitions.

\begin{Definition}\label{def:1} Let $A(\mu,\nu,N)$ denote the number of
partitions of $N$ into distinct non-multiples of 6 of which $ \mu $
are congruent to 1 or 2 mod~6
and $\nu $  are congruent to 4 or 5 mod~6.
\end{Definition}

Clearly, we have
\begin{equation} \label{eq:GF}
\sum_{\mu,\nu,N \geq 0}A(\mu,\nu,N)a^\mu b^\nu q^N
 =
\prod^\infty_{n=0}(1+aq^{6n+1})(1+aq^{6n+2})(1+bq^{6n+4})(1+bq^{6n+5}).
\end{equation}

\begin{Definition}\label{def:2}
Let $B(\mu,\nu,N)$ denote the number of partitions
$ \lambda = b_1+\cdots+b_s $ of N satisfying the following conditions:

\begin{itemize}
\item [(i)]Only multiples of 6 may be repeated.
\item [(ii)]$ b_i-b_{i+2} \geq 6 $ with strict inequality if $b_i$ is
a multiple of 6.
\item [(iii)]The multiplicities $f_i$, $1\le i\le s$, satisfy
\begin{align*}
  f_{6j+3}&=0 \quad \text{for all}\quad  j \geq 0 ,\\
       f_{6j+2}+f_{6j+4} &\leq 1 \quad \text{for all}\quad j \geq 0 ,\\
       f_{6j+5}+f_{6j+7} &\leq 1 \quad \text{for all}\quad j \geq 0 ,\\
       f_{6j-1}+f_{6j}+f_{6j+6}+f_{6j+7} &\leq 3\quad \text{for all}\quad
j \geq 1 .
\end{align*}
\item [(iv)]There are $\mu$ parts of the partition $ \equiv $ 0, 1 or 2
(mod 6).
\item [(v)] There are $\nu$ parts of the partition $ \equiv $ 0, 4 or 5
(mod 6).
\end{itemize}
\end{Definition}

The following theorem is the announced refinement of Theorem~\ref{thm:2}.

\begin{Theorem}\label{thm:3}
For each $\mu,\nu,N \geq 0 $ we have
\begin{equation} \label{eq:part}
 A(\mu,\nu,N) = B(\mu,\nu,N).
\end{equation}
\end{Theorem}

It is obvious that Theorem~\ref{thm:2} follows immediately from
Theorem~\ref{thm:3} by summing both sides of \eqref{eq:part} over all
$\mu$ and $\nu$. The proof of Theorem~\ref{thm:3} is given in the next
section.

\section{Proof of Theorem~\ref{thm:3}} \label{sec:3}

We first observe that for any partition which satisfies
(i)--(iii) of Definition~\ref{def:2} there are exactly
16 possibilities (numbered 0--15 in Table~\ref{tab:1})
for the subset of summands of the partition that lie in the interval
$[6i+1, 6i+6]$.

\begin{table}[h]
$$
\begin{array}{rl}
0:\quad& \emptyset : \text{empty} \\
 1: \quad&  6i+1 \\
 2:\quad&  6i+2 \\
 3:\quad&  6i+2,\ 6i+1 \\
 4: \quad& 6i+4 \\
 5: \quad& 6i+4,\ 6i+1 \\
 6:\quad&  6i+5 \\
 7: \quad& 6i+5,\ 6i+1 \\
 8:\quad&  6i+5,\ 6i+2 \\
 9: \quad& 6i+5,\ 6i+4 \\
 10:\quad&  6i+6 \\
 11:\quad&  6i+6,\ 6i+1 \\
 12:\quad&  6i+6,\ 6i+2 \\
 13: \quad& 6i+6,\ 6i+4 \\
 14:\quad&  6i+6,\ 6i+5 \\
 15: \quad& 6i+6,\ 6i+6
\end{array}
$$
\caption{}
\label{tab:1}
\end{table}

We now refine the partitions from Definition~\ref{def:2} further,
using the classification given in Table~\ref{tab:1}.

\begin{Definition}\label{def:3}
Let $S_n(j,a,b,q)$ denote the generating function
$$
\sum _{} ^{}B_n(\mu,\nu,N)a^\mu b^\nu q^N,$$
where $B_n(\mu,\nu,N)$ is the number of all
partitions considered in Definition~\ref{def:2}, which in
addition satisfy the two conditions
\begin{itemize}
\item [(vi)] all parts are $\leq 6n+6$, and
\item[(vii)] the subset of summands
that lie in the interval $[6n+1, 6n+6]$ must have number $\leq j$
in Table~\ref{tab:1}.
\end{itemize}
When $n = -1$, we define $ S_{-1}(j,a,b,q) = 1 $ and for
$ n < -1,$ we define $ S_n(j,a,b,q) = 0. $
\end{Definition}

For example,
$$
S_0(9,a,b,q) = 1+aq+aq^2+a^2q^3+bq^4+abq^5+bq^5+abq^6+abq^7 $$
and
\begin{multline*}
 S_0(15,a,b,q) =
1+aq+aq^2+a^2q^3+bq^4+abq^5+bq^5+2abq^6+a^2bq^7+abq^7\\
+a^2bq^8+b^2q^9
 +ab^2q^{10}+ab^2q^{11}+a^2b^2q^{12}.
\end{multline*}
 It is easy to verify that
$$ S_0(15,a,b,q) =
(1+aq)(1+aq^2)(1+bq^4)(1+bq^5) .$$

For convenience, we write
$S_n(j)$ for $S_n(j,a,b,q)$. Along the lines of \cite{R3}, we obtain
the following recurrence relations for $S_n(j)$:
\begin{align} \label{eq:16}
S_n(0)&=S_{n-1}(15),\\
\label{eq:17} S_n(1) &=
S_n(0)+aq^{6n+1}[S_{n-1}(11)-S_{n-1}(9)+S_{n-1}(5)]\\
\notag
&\kern6cm
-a^3b^3q^{24n-12}S_{n-3}(9),\\
\label{eq:18} S_n(2) &=
S_n(1)+aq^{6n+2}[S_{n-1}(12)-S_{n-1}(9)+S_{n-1}(8)], \\
\label{eq:19} S_n(3) &= S_n(2)+a^2q^{12n+3}S_{n-1}(3),\\
\label{eq:20} S_n(4) &= S_n(3)+bq^{6n+4}S_{n-1}(13), \\
\label{eq:21} S_n(5) &= S_n(4)+abq^{12n+5}S_{n-1}(5), \\
\label{eq:22} S_n(6) &= S_n(5)+bq^{6n+5}S_{n-1}(14), \\
\label{eq:23} S_n(7) &= S_n(6)+abq^{12n+6}S_{n-1}(5), \\
\label{eq:24} S_n(8) &= S_n(7)+abq^{12n+7}S_{n-1}(8), \\
\label{eq:25} S_n(9) &= S_n(8)+b^2q^{12n+9}S_{n-1}(9), \\
\label{eq:26} S_n(10)&= S_n(9)+abq^{6n+6}S_{n-1}(14), \\
\label{eq:27} S_n(11) &= S_n(10)+a^2bq^{12n+7}S_{n-1}(5), \\
\label{eq:28} S_n(12) &= S_n(11)+a^2bq^{12n+8}S_{n-1}(8), \\
\label{eq:29} S_n(13) &= S_n(12)+ab^2q^{12n+10}S_{n-1}(9),\\
\label{eq:30} S_n(14) &= S_n(13)+ab^2q^{12n+11}S_{n-1}(9), \\
\label{eq:31} S_n(15) &= S_n(14)+a^2b^2q^{12n+12}S_{n-1}(9).
\end{align}

\medskip
We now define two linear combinations of the $S_n(9)$'s and the $S_n(15)$'s,
\begin{align} \label{eq:32}
  J(n): =
S_n&(9)-(1-q^{6n})(1+aq^{6n+1}+aq^{6n+2}+bq^{6n+4}+bq^{6n+5})S_{n-1}(15)\\
\notag
&-q^{6n}(1+aq^{6n+1}+aq^{6n+2}+a^2q^{6n+3}+bq^{6n+4}+bq^{6n+5}+abq^{6n+5}\\
\notag
&\kern5cm
+abq^{6n+6}+abq^{6n+7}+b^2q^{6n+9})S_{n-1}(9)\\
\notag
&+(1-q^{6n})abq^{18n-3}
(a^2+abq^2+abq^3+abq^4+a^2bq^4+a^2bq^5+b^2q^6\\
\notag
&\kern5cm+ab^2q^7+ab^2q^8)S_{n-2}(9)\\
\notag & +a^3b^3q^{24n-12}(1-q^{6n})(1-q^{6n-6})S_{n-3}(9),
\end{align}
and
\begin{align} \label{eq:33}
K(n) : =  S_n(9) - &S_n(15) + abq^{6n+6}(1 - q^{6n})S_{n-1}(15)\\
\notag
& + abq^{12n+6} (1+aq+aq^2+bq^4+bq^5+abq^6) S_{n-1}(9)\\
\notag
&-a^3b^3q^{18n+6}(1-q^{6n})S_{n-2}(9).
\end{align}

Along lines similar to those in \cite{R3}, we are able to
obtain a recurrence for $S_n(9)$ (see Lemma~\ref{lem:2}).

\begin{Lemma}\label{lem:1}
For $n \geq  0$, $J(n)  =  0  =  K(n)$.
\end{Lemma}

\begin{proof}[Sketch of proof]
We prove the lemma by using the 
identities \eqref{eq:16}--\eqref{eq:31}. The 14
sequences $S_n(j)$, where $j$ is different from $9$ and $15$, 
can be expressed as linear combinations of the $S_m(9)$'s and
the $S_m(15)$'s in the following way: from \eqref{eq:31}, 
we find that $S_n(14)$ is given by
$$S_n(14)=S_n(15)-a^2b^2q^{12n+12}S_{n-1}(9).$$
Using the above equation in \eqref{eq:30}, $S_n(13)$ becomes such a linear
combination. Similarly for $S_n(12)$ if we use \eqref{eq:29}. 
Equation \eqref{eq:25}
yields such a linear combination for $S_n(8).$
Subsequently, \eqref{eq:28} yields a linear combination for 
$S_n(11)$, and \eqref{eq:26} yields a linear combination for $S_n(10).$
Replacing $n$ by $n+1$ in \eqref{eq:27}, we get 
$$S_n(5) =
a^{-2}b^{-1}q^{-12n-19}[S_{n+1}(11)-S_{n+1}(10)],$$ 
which in turn yields an
expression for $S_n(5)$ in terms of the $S_m(9)$'s and the $S_m(15)$'s. 
Equations~\eqref{eq:24}, \eqref{eq:23},  \eqref{eq:21},
 \eqref{eq:20},  \eqref{eq:19} and  \eqref{eq:18}
yield respectively linear combinations 
in terms of the $S_m(9)$'s and the $S_m(15)$'s for
$S_n(7)$, $S_n(6)$, $S_n(4)$, $S_n(3)$, $S_n(2)$ and $S_n(1)$. 
Finally, we know already from
 \eqref{eq:16} that $S_n(0) = S_{n-1}(15)$.

We are now in the position to prove $K(n) = 0.$
Let us consider  \eqref{eq:22}, that is 
$$S_n(6) = S_n(5)+bq^{6n+5}S_{n-10}(14).$$
Substituting the expression in terms of the $S_m(9)$'s
and $S_m(15)$'s obtained earlier for $S_n(5)$ and
the respective one for $S_{n-1}(14)$ in the equation above, 
we get a certain identity, (A) say.

On the other hand, from \eqref{eq:23}, we have 
$$S_n(6) = S_n(7)-abq^{12n+6}S_{n-1}(5).$$
Substituting the expression in terms of the $S_m(9)$'s
and $S_m(15)$'s obtained earlier for 
$S_n(7)$ and the respective one for $S_{n-1}(5)$, 
we obtain another identity, (B) say.
It can now be verified that (A)$-$(B), when multiplied by 
$a^2bq^{12n+19}$, is exactly the equation $K(n+1) = 0.$ 

Now we prove $J(n) = 0.$ 
Substituting the expressions obtained earlier
for $S_n(1)$, $S_n(0)$, $S_{n-1}(11)$
and $S_{n-1}(5)$ into \eqref{eq:17}, we obtain
$$ 0 =
a^2bq^{12n+19}J(n)-K(n+1)+aq^{6n+13}(1+aq^{6n+2}+bq^{6n+4}+bq^{6n+5})K(n).$$
Since $K(n) = 0$ for all $ n \geq 0,$ we conclude that $J(n) = 0$
for all $n\ge0$. This proves the lemma.
\end{proof}

\begin{Lemma}\label{lem:2}
For  $n\geq 0$,
\begin{align} \label{eq:Id}
 (1+a&q^{6n-5}+aq^{6n-4}+bq^{6n-2}+bq^{6n-1})S_n(9)\\
\notag
& = p_1(n,a,b,q)S_{n-1}(9)+(1-q^{6n})p_2(n,a,b,q)S_{n-2}(9)\\
\notag
&\kern1cm
+p_3(n,a,b,q)(1-q^{6n})
(1-q^{6n-6})S_{n-3}(9)\\
\notag
&\kern1cm
+a^4b^4q^{30n-36}(1-q^{6n})(1-q^{6n-6})(1-q^{6n-12})\\
\notag
&\kern2cm
\times[(1+aq^{6n+1}+aq^{6n+2}+bq^{6n+4}+bq^{6n+5})S_{n-4}(9)],
\end{align}
where
\begin{align} \label{eq:34}
 p_1(n,a,b,q) ={}&
1+aq^{6n-5}+aq^{6n-4}+bq^{6n-2}+bq^{6n-1}+abq^{6n}+2abq^{12n-1}+3abq^{12n}\\
\notag
&
+aq^{6n+2}+2a^2q^{12n-3}+a^2q^{12n-2}+2abq^{12n+1}+a^2bq^{12n+2}+bq^{6n+4}\\
\notag
&
+b^2q^{12n+2}+2b^2q^{12n+3}+ab^2q^{12n+4}+bq^{6n+5}+b^2q^{12n+4}+ab^2q^{12n+5}\\
\notag
&
+a^2q^{12n+3}+2a^2bq^{18n+1}+2a^2bq^{18n+2}+abq^{12n+5}+a^2bq^{18n}\\
\notag
& +ab^2q^{18n+3}+a^2q^{12n-4}+2ab^2q^{18n+4}+abq^{12n+6}+2ab^2q^{18n+5}\\
\notag
&+abq^{12n+7}+a^2bq^{18n+3}+ab^2q^{18n+6}+b^2q^{12n+9}+a^3q^{18n-2}\\
\notag
&+a^3q^{18n-1}+b^3q^{18n+7}+b^3q^{18n+8}+a^2bq^{12n+1}+aq^{6n+1},
\end{align}
\begin{align}\label{eq:35}
 p_2(n,a,b,q) ={}&
a^2bq^{12n-5}+a^2bq^{12n-4}+ab^2q^{12n-2}+ab^2q^{12n-1}+a^2b^2q^{12n}+a^3bq^{18n-4}\\
\notag
&+a^3bq^{18n-3}+a^3bq^{18n-2}+3a^2b^2q^{18n}+a^2b^2q^{18n-1}+ab^3q^{18n+2}\\
\notag
&+ab^3q^{18n+3}+a^2b^2q^{18n+1}+ab^3q^{18n+4}+a^3bq^{18n-10}+a^2b^2q^{18n-7}
\\
\notag
&+3a^2b^2q^{18n-6}+a^3b^2q^{18n-5}+a^4bq^{24n-9}+a^4bq^{24n-8}+a^4bq^{24n-7}\\
\notag
&+3a^3b^2q^{24n-5}+a^3b^2q^{24n-6}+3a^2b^3q^{24n-2}+3a^3b^2q^{24n-4}+3a^2b^3q^{24n-1}\\
\notag
&+a^3bq^{18n-9}+a^3bq^{18n-8}+a^2b^2q^{18n-5}+a^3b^2q^{18n-4}+a^2b^3q^{24n}\\
\notag
&+ab^3q^{18n-4}+ab^4q^{24n}+a^3b^2q^{24n-3}+ab^4q^{24n+2}+ab^3q^{18n-2}\\
\notag
&+ab^4q^{24n+1}+ab^4q^{24n+3}+a^4bq^{24n-6}+ab^3q^{18n-3}+a^2b^3q^{18n-2}\\
\notag
&+a^2b^3q^{18n-1}+a^2b^3q^{24n-3},
\end{align}
and
\begin{align} \label{eq:36}
 p_3(n,a,b,q) ={}&
-a^3b^3q^{24n-12}-a^3b^3q^{18n-12}-a^4b^3q^{24n-11}-a^4b^3q^{24n-10}-a^3b^4q^{24n-8}\\
\notag
&-a^3b^4q^{24n-7}+2a^4b^3q^{30n-17}+2a^4b^3q^{30n-16}+2a^3b^4q^{30n-14}\\
\notag
&+2a^3b^4q^{30n-13}+a^4b^2q^{24n-21}+a^5b^2q^{30n-20}+a^3b^3q^{24n-19}\\
\notag
&+a^5b^2q^{30n-19}+a^4b^3q^{30n-18}+a^3b^4q^{30n-15}+a^3b^3q^{24n-18}+a^3b^3q^{24n-17}\\
\notag
&+a^4b^3q^{30n-15}+a^3b^4q^{30n-12}+a^2b^5q^{30n-11}+a^2b^5q^{30n-10}+a^2b^4q^{24n-15}.
\end{align}
\end{Lemma}

\begin{proof}[Sketch of proof]
Using the identity $J(n) = 0$, we find that $S_{n-1}(15)$ is a linear
combination of $S_n(9)$, $S_{n-1}(9)$, $S_{n-2}(9)$ and $S_{n-3}(9).$
By Lemma~5, we have $K(n) = 0.$ We substitute that linear combination
for $S_{n-1}(15)$ and the corresponding one for
$S_{n}(15)$ in \eqref{eq:33}. 
After some simplification, and after replacing $n$ by $n-1$,
we arrive exactly at \eqref{eq:Id}.
\end{proof}

We are now able to prove a recurrence for $S_n(15)$.

\begin{Lemma}\label{lem:3}
For $n\geq 0$, we have
\begin{align} \label{eq:37}
 (1+a&q^{6n-5}+aq^{6n-4}+bq^{6n-2}+bq^{6n-1})S_n(15)\\
\notag
&=
p_1(n-1,aq^6,bq^6,q)S_{n-1}(15)+(1-q^{6n-6})p_2(n-1,aq^6,bq^6,q)S_{n-2}(15)\\
\notag
&\kern1cm
+(1-q^{6n-6})(1-q^{6n-12})p_3(n-1,aq^6,bq^6,q)S_{n-3}(15)\\
\notag
&\kern1cm
+a^4b^4q^{30n-18}
(1+aq^{6n+1}+aq^{6n+2}+bq^{6n+4}+bq^{6n+5})\\
\notag
&\kern2cm\times(1-q^{6n-6})(1-q^{6n-12})(1-q^{6n-18})S_{n-4}(15).
\end{align}
\end{Lemma}

\begin{proof}
Since $J(n) = 0$, we can express $S_{n}(15)$ in terms of the $S_{n}(9)'s.$
Substituting these expressions in \eqref{eq:37},
we get an equation involving $S_{n+1}(9),S_{n}(9),\dots,
S_{n-6}(9).$ In that equation we apply Lemma~\ref{lem:2}
to $S_{n+1}(9).$ In the result thus obtained,
we apply Lemma~\ref{lem:2} to $S_{n}(9).$ In the subsequent result
obtained, we apply Lemma~\ref{lem:2} to $S_{n-1}(9).$
In the result obtained, we again apply Lemma~\ref{lem:2}, this time
to $S_{n-2}(9)$. The result is zero. All these calculations
have been performed using {\sl Mathematica}.
\end{proof}

\begin{Lemma}\label{lem:4}
For $n\geq0$,
$$ S_n(15,a,b,q) =
(1+aq)(1+aq^2)(1+bq^4)(1+bq^5)S_{n-1}(9,aq^6,bq^6,q). $$
\end{Lemma}

\begin{proof}
Comparing Lemma~\ref{lem:3} and Lemma~\ref{lem:2}, we find that
both sides of Lemma~\ref{lem:4} satisfy the same fourth order
recurrence valid for $ n \geq 1. $ Hence we have only to verify
Lemma~\ref{lem:4} for $n = 0,1,2,3$. This is a routine
verification, and can therefore be left to the reader.
\end{proof}

\begin{proof}[Proof of Theorem~\ref{thm:3}]
For $  0 \leq j \leq 15, $ we have
$$ \lim_{n\to\infty}S_{n}(j,a,b,q)
 = \sum_{\mu,\nu,N \geq 0} B(\mu,\nu,N)a^{\mu}b^{\nu}q^N \equiv S(a,b,q).$$
 Letting $ n\rightarrow \infty $ in Lemma~\ref{lem:4}, we find that
 $$ S(a,b,q) =(1+aq)(1+aq^2)(1+bq^4)(1+bq^5)S(aq^6,bq^6,q) .$$
 Iterating the above equation, we obtain
\begin{align*}
 S(a,b,q) &= \prod_{n=0} ^\infty
(1+aq^{6n+1})(1+aq^{6n+2})(1+bq^{6n+4})(1+bq^{6n+5}) \\
& = \sum_{\mu,\nu,N \geq 0} A(\mu,\nu,N)a^{\mu}b^{\nu}q^N ,
\end{align*}
the latter equality being due to \eqref{eq:GF}.
 Thus we get
 $$  \sum_{\mu,\nu,N \geq 0} B(\mu,\nu,N)a^{\mu}b^{\nu}q^N
=  \sum_{\mu,\nu,N \geq 0} A(\mu,\nu,N)a^{\mu}b^{\nu}q^N ,$$
 which, upon comparison of coefficients of $a^{\mu}b^{\nu}q^N $, implies
  $$ A(\mu,\nu,N) = B(\mu,\nu,N) $$
  for all non-negative $\mu,\nu$ and $N$. This is exactly the claim in
Theorem~\ref{thm:3}.
\end{proof}

\end{document}